# ADAPTIVE INDEPENDENT METROPOLIS–HASTINGS[1]

By Lars Holden, Ragnar Hauge and Marit Holden

*Norwegian Computing Center*

We propose an adaptive independent Metropolis–Hastings algorithm with the ability to learn from all previous proposals in the chain except the current location. It is an extension of the independent Metropolis–Hastings algorithm. Convergence is proved provided a strong Doeblin condition is satisfied, which essentially requires that all the proposal functions have uniformly heavier tails than the stationary distribution. The proof also holds if proposals depending on the current state are used intermittently, provided the information from these iterations is not used for adaption. The algorithm gives samples from the exact distribution within a finite number of iterations with probability arbitrarily close to 1. The algorithm is particularly useful when a large number of samples from the same distribution is necessary, like in Bayesian estimation, and in CPU intensive applications like, for example, in inverse problems and optimization.

**1. Introduction.** Assume we want to sample from a distribution $\pi$ or make estimates based on $\pi$, but direct samples from $\pi$ are not obtainable. Rejection sampling importance sampling, and sampling importance resampling (SIR) are techniques for generating such samples and estimates by proposing samples from a different distribution $q$. Another alternative is to use a Metropolis–Hastings algorithm with a proposal function that is independent of the present position. This approach, which we call independent Metropolis–Hastings, is also known as independent Markov chain in Tierney (1994), or Metropolized independent sampling in Liu (1996), or independence sampler in Roberts and Rosenthal (1998). The efficiency of these methods depends on the proposal distribution $q$ being close to $\pi$. If this is not practically possible, other Markov chains such as Metropolis–Hastings based on local moves around the present state may be better; see Meyn and Tweedie (1993), Gilks et al. (1996) and Geyer (1992).

Received August 2006; revised February 2008.

[1]Supported by Research Council of Norway.

*AMS 2000 subject classifications.* Primary 65C05; secondary 65C40.

*Key words and phrases.* Adaption, Metropolis–Hastings, Markov chain Monte Carlo, inverse problems.

In all the alternatives described above, knowledge about the stationary distribution is gained as the number of iterations increases. This knowledge may be used to adapt the proposal distribution in order to improve the convergence of the chain. The bound on the convergence improves and the bound on the correlation between subsequent states in the chain decreases when the proposal distribution better approximates the stationary distribution; see the proposition in Holden (1998b).

Adaptive Markov chains is an active research area; see, for example; Atchade and Rosenthal (2003) and the references therein. The paper Roberts and Rosenthal (2005) has several general results for adaptive Markov chains, but they focus on diminishing adaption that is not necessarily satisfied by the method proposed in this paper. Erland (2003) gives an overview of adaptive algorithms. He divides adaptive algorithms into three groups: adaptive strategies within ordinary MCMC, algorithms where the adaption diminishes and algorithms with regeneration. There are several papers proposing very different algorithms within each group. The possibilities for adaption within ordinary MCMC is limited; see Tjelmeland and Hegstad (2001). Regeneration tends to happen so seldom in higher dimension that this limits the applicability of these algorithms; see Gilks et al. (1998). Algorithms where the adaption diminishes are probably the most promising; see for example, Haario et al. (2001). However, also these algorithms have restrictions that make them difficult to use.

The algorithm presented in this paper is within the diminishing subgroup even though the adaption does not need to diminish. It is less technical than the other algorithms, and it is easy to describe and use in practice. Also, it has a proposal function that depends on all earlier proposed states except the current, and a Metropolis–Hastings-like acceptance step. There is a large flexibility in the choice of proposal function and it is a challenge to find a proposal function that is able to use the data efficiently in order to obtain satisfactory convergence. This is discussed in the paper and four different examples are given with both parametric and nonparametric alternatives. We call this algorithm the adaptive independent Metropolis–Hastings algorithm. This is a special case of the adaptive chain described in Holden (1998a). Surprisingly, the limiting density is invariant with an independent proposal function and hence this gives much better convergence properties than in the general case. The algorithm is particularly designed for examples where it is expensive to calculate the limiting density $\pi$. In such examples it is most likely cost effective to use a lot of data in an adaptive scheme.

Gåsemyr (2003) describes another adaptive algorithm based on the independent Metropolis–Hastings algorithm. The convergence rate of both that algorithm and our algorithm depends on the ratio $\pi/q_i$, where $q_i$ denotes the proposal distribution in iteration $i$. Both algorithms do also give exact samples in a finite number of iterations. The adaption schemes seem different, as the algorithm in Gåsemyr (2003) adapts on previously accepted states, whereas our algorithm adapts on proposed states. However, this difference is not as essential since the algorithm presented in Gåsemyr (2003) may easily use proposed states instead. From all the proposed states, it is possible to generate a new independent chain with the same properties as the accepted chain. The main difference is that the algorithm described in Gåsemyr (2003) requires that the supremum of the ratio $f/q_i$ is known, where $f = c\pi$ and $c$ is an unknown constant. This supremum is used in the algorithm, and although the algorithm will converge without it, the convergence will be weaker. It adapts using groups of previous states, and the adaption stops after a finite number of iterations. The algorithm presented in this paper only requires that the supremum of $\pi/q_i$ is finite, it may use all previous states in the adaption and the adaption does not stop. It may also include iterations with proposal functions that depend on the present position, but the information from these iterations cannot be used for adaption. Andrieu and Moulines (2003) also discuss independent Metropolis–Hastings and our Theorem 1 is of interest for the research area they focus on.

The independent Metropolis–Hastings algorithm is an efficient sampling algorithm only if $q_i$ is reasonably close to $\pi$. However, it allows great freedom of adaption, as shown both here and in Gåsemyr (2003). Adaption should also make the proposals resemble $\pi$, and hence the adaptive version will be an efficient sampler. We will prove a bound on the convergence that depends on the supremum of $\pi/q_i$. An attractive property of independent proposals is their ability to make large jumps, and if this can be done while keeping the acceptance rate high, the autocorrelation of the chain will decrease rapidly.

**2. Definition of the adaptive independent Metropolis–Hastings algorithm.**
The goal for the algorithm is to sample from a distribution $\pi$, which is known except for a normalizing constant. The algorithm resembles the traditional Metropolis–Hastings algorithm (MH). A chain $x$ is generated by drawing proposals $z_i$ from a proposal distribution $q_i(z_i|x_{i-1}, \tilde{y}^{i-1})$ and either accepting them, setting $x_i = z_i$, or rejecting them, setting $x_i = x_{i-1}$. The proposal history vector $\tilde{y}^i$ is defined as follows: If the proposal is independent of the current state, the history vector is extended by including $z_i$ if $z_i$ was rejected, and $x_{i-1}$ otherwise. On the other hand, if the proposal is dependent on the current state, the history vector is kept unmodified, that is, $\tilde{y}^i = \tilde{y}^{i-1}$. The difference between traditional MH and adaptive independent MH is only that the proposal function $q_i$ may depend on the history vector. In adaptive independent MH $q_i$ may depend on all states where the function $f = c\pi$ has been evaluated, except the current state of the chain, and these values can be used to make $q_i$ a better approximation of $\pi$. Also when doing local steps conditioned on $x_i$, all information in $\tilde{y}^i$ can be used. The limitation is that information gained from doing the local steps cannot be used, so these iterations do not improve the proposal.

The full algorithm is given as:

1. Set $\tilde{y}^0 = \varnothing$.
2. Generate an initial state $x_0 \in \Omega$ from the density $p_0$.
3. For $i = 1, \ldots, n$:
   (a) Generate a state $z_i$ from the density $q_i(z_i | x_{i-1}, \tilde{y}^{i-1})$.
   (b) Calculate $\alpha_i(z_i, x_{i-1}, \tilde{y}^{i-1}) = \min\{1, \frac{\pi(z_i) q_i(x_{i-1} | z_i, \tilde{y}^{i-1})}{\pi(x_{i-1}) q_i(z_i | x_{i-1}, \tilde{y}^{i-1})}\}$.
   (c) Set

   $$x_i = \begin{cases} z_i, & \text{with probability } \alpha_i(z_i, x_{i-1}, \tilde{y}^{i-1}), \\ x_{i-1}, & \text{with probability } 1 - \alpha_i(z_i, x_{i-1}, \tilde{y}^{i-1}). \end{cases}$$

   (d) Set $\tilde{y}^i = \tilde{y}^{i-1}$ if $q_i$ depended on $x_{i-1}$. Otherwise,

   $$\tilde{y}^i = \tilde{y}^{i-1} \text{ appended with } \begin{cases} z_i, & \text{if } z_i \text{ was rejected}, \\ x_{i-1}, & \text{if } z_i \text{ was accepted.} \end{cases}$$

A classical random-walk MH algorithm where only states close to the present state in the chain are proposed, may stay close to a local optimum for a large number of iterations such that the user believes that the chain has converged. This may also be a problem with the presented algorithm, but if the proposal function is used properly, the probability for staying close to a local optimum is less with this algorithm than with classical MH since more information may be included in the proposal. However, if the proposal function adapts too fast, the presented algorithm may be even worse than classical MH.

The presented algorithm gives a very large flexibility in choosing the proposal function $q_i$. Creating a good proposal function may be difficult. The choice depends on the problem we want to solve, the function $\pi$, the time used for an evaluation of $\pi$, how to let $q_i$ approximate $\pi$, and the CPU and software resources available. It may also be a challenge to use all the data $\tilde{y}^i$ effectively, in particular in higher dimensions. Theorem 2 gives a bound on the convergence based on how well $q_i$ approximates $\pi$, similar to the convergence bounds for the Metropolis–Hastings algorithm given in Holden (1998b). High convergence rate also implies that the correlation between state $x_i$ and $x_{i+k}$ decreases quickly as $k$ increases. This may be as important as convergence if the chain is used to compute averages such as

$$\frac{1}{N} \sum_{i=1}^{N} F(x_i)$$

for some function $F$.

Finding a suitable representation of $q_i$ that is able to use all the data $\tilde{y}^i$ effectively may be a challenge. For a parametric $q_i$ this is done by estimating the parameters from the knowledge about $\pi$ in the history states. This is illustrated in Example 4. However, $\pi$ may be too complicated for a simple parametric description. In these cases, nonparametric approaches can be used, as is done in three of the examples. It is always important that the proposal function does not adapt too fast and has heavy tails. This becomes clear when we see the convergence properties and the examples. Another problem is that the history vector grows large, and using all the information there for an update may become time consuming. Again, the method used in the examples shows how this information may be reduced. In cases where the evaluation of $f = c\pi$ is expensive, the overhead in computing the proposal even with a full history vector may be insignificant. Kriging or other spatial statistical methods are possible, but these alternatives may be complicated if the number data points become large. In a large class of problems [see, e.g., Park and Jeon (2002) and Kleijnen and Sargent (2000)],

$$f(x) = \sum_{j=1}^{m} \omega_j (g_j(x) - d_j)^2,$$

where $(g_1, g_2, \ldots, g_m)$ is the vector-valued response of a simulation model, $\omega_j$ are weights and $d_j$ are the data. This is an inverse problem where it is necessary to find the parameter values $x$ with uncertainty that give an observed response. For these problems we may use the evaluations of $g_j(z_i)$ to find an approximation to the functions $g_j(x)$ parametrically or nonparametrically. Example 4 illustrates this approach.

**3. Convergence.** Let $\Omega^1 = \Omega \subset \mathbb{R}^n$ be a Borel-measurable state space or, alternatively, let $\Omega$ be a discrete state space, and $\Omega^{i+1} = \Omega \times \Omega^i$, and $\tilde{x}^i = (x_1, \ldots, x_i) \in \Omega^i$. Let $\mu(\tilde{x}^i)$ be the product measure on $\Omega^i$ of a $\sigma$ finite measure $\mu(x)$. Further, let the density $\pi$ and the proposal functions $q_i$ be integrable with respect to $\mu$ including point mass distribution. In the notation we neglect that $\tilde{y}^i$ may have dimension less than $i$ in order to simplify the notation. This is the case if some of the proposal functions $q_i$ use local proposals.

The convergence rate depends on the constant $a_i(\tilde{y}^{i-1})$ in the strong Doeblin condition: Let $a_i(\tilde{y}^{i-1}) \in [0, 1]$ satisfy

(1) $\quad q_i(z|x, \tilde{y}^{i-1}) \geq a_i(\tilde{y}^{i-1})\pi(z) \qquad$ for all $(z, x, \tilde{y}^{i-1}) \in \Omega^{i+1}$ and all $i > 1$.

The Doeblin condition essentially requires that all the proposal distribution has uniformly heavier tails than the target distribution. This condition is always satisfied for $a_i(\tilde{y}^{i-1}) = 0$. Theorem 2 below is valid also in this case, but we only prove convergence if the expected value satisfies $E(\prod_{i=1}^{\infty}(1 - a_i(\tilde{y}^{i-1}))) = 0$. However, it may be useful to allow that $a_i(\tilde{y}^{i-1}) > 0$ for some values of $i$. It is possible to obtain $a_i(\tilde{y}^{i-1}) > 0$ by making the tails in the proposal function sufficiently heavy. This will usually be satisfied in a Bayesian approach using a prior with heavy tails. Alternatively, we may replace $q_i(z|x, \tilde{y}^{i-1})$ by $(1 - \varepsilon)q_i(z|x, \tilde{y}^{i-1}) + \varepsilon g(z)$ where $g(z)$ is a density with extremely heavy tails and $\varepsilon$ is small. The assumptions in the theorem are satisfied and we only "waste" at most $\varepsilon$ of the proposals. However, a small $\varepsilon$ may give slow convergence in the tails. The Doeblin condition is natural for an independent sampler. If it is not satisfied for some states in $\Omega$, the algorithm will tend to undersample these states. This is not crucial if the probability mass of these states is low and further inference does not depend on tail behavior. In one of our examples the stationary distribution has heavier tails than the proposal distribution, but convergence is still achieved with the measure we use.

Let $p_i(x_i)$ be the distribution for $x_i$ after $i$ iterations and $\hat{p}_i(x_i|\tilde{y}^i)$ the conditional distribution for $x_i$ given the history vector. The following theorem is crucial in the understanding of the algorithm, and is the basis for the convergence result.

THEOREM 1. *The limiting density $\pi$ conditioned on the history is invariant for the adaptive independent Metropolis–Hastings algorithm, that is, $\hat{p}_i(x_i|\tilde{y}^i) = \pi(x_i)$ implies $\hat{p}_{i+1}(x_{i+1}|\tilde{y}^{i+1}) = \pi(x_{i+1})$.*

It is possible to combine this theorem with Theorem 6 in Roberts and Rosenthal (2005) to get the following result: If we assume that the proposal function for all histories is uniform ergodic and the adaption diminishes, then the algorithm converges in the total variance (TV) norm

$$\|p_i - \pi\|_{\text{TV}} = \int_{\Omega} |p_i(x) - \pi(x)| \, d\mu(x).$$

In this paper we will instead assume the strong Doeblin condition and obtain geometric convergence.

THEOREM 2. *Assume the adaptive independent Metropolis–Hastings algorithm satisfies (1). Then*

$$(2) \qquad \|p_i - \pi\|_{\text{TV}} \leq 2 \int_{\Omega^i} \tilde{p}_i(\tilde{y}^i) \prod_{j=1}^{i} (1 - a_j(\tilde{y}^{j-1})) \, d\mu(\tilde{y}^{i-1})$$

$$(3) \qquad = 2E\left(\prod_{j=1}^{i} (1 - a_j(\tilde{y}^{j-1}))\right).$$

*If $a_j(\tilde{y}^{j-1}) \geq a_j$ for all $j$ and $\tilde{y}^{j-1} \in \Omega^{j-1}$, then*

$$(4) \qquad \|p_i - \pi\|_{\mathrm{TV}} \leq 2 \prod_{j=1}^{i}(1 - a_j).$$

*The algorithm converges if this product goes to zero when $i \to \infty$. If $a_j > a > 0$ infinitely often, the algorithm samples from the target distribution within a finite number of samples with a probability arbitrarily close to 1.*

The proofs are given in the Appendix. Theorem 2 says that convergence is geometric as long as a strong Doeblin condition is satisfied for all possible histories with $a_j(\tilde{y}^{j-1}) \geq a > 0$ for all $j$ and $\tilde{y}^{j-1} \in \Omega^{j-1}$. In each iteration the chain jumps to the limiting density $\pi$ with probability $a_j(\tilde{y}^{j-1})$. Then the chain remains in this density according to Theorem 1 since $\pi$ is invariant for the adaptive independent Metropolis–Hastings algorithm. The probability for not sampling from the limiting density after $j$ iterations is $E(\prod_{j=1}^{i}(1 - a_j(\tilde{y}^{j-1})))$.

If the adaption succeeds in generating better proposal distributions, $a_i(\tilde{y}^{i-1})$ will increase as $i$ increases, and the convergence will be accelerating. This also means that the number of iterations needed to generate a set of independent samples decreases.

**4. Examples.** This section gives four examples. All examples are schematic examples where we compare different Metropolis–Hastings algorithms: independent sampler, random-walk proposing small jumps and adaptive independent sampler where the proposal function is higher close to the modes that are identified so far. The first example is a quantitative comparison showing that the adaptive independent sampler converges faster, identifies modes better and jumps more often between different modes than the other algorithms. In Examples 2 and 3 the same adaptive independent algorithm is used showing how flexibly the algorithm may be used in different cases. The three first examples are in one or two dimensions. We do not believe that the dimension is critical. What is important, though, is how sharp the modes are and the number of modes. But in high dimension one may prefer to change only a few variables in each iteration if this reduces the number of calculations per iteration and not because this is an efficient Metropolis–Hastings algorithm per iteration. The final example illustrates how to combine the proposed algorithm with an external simulation model that is assumed to be very demanding to evaluate. This is a typical problem in a large number of applications.

EXAMPLE 1. Let $\pi$ be the function

$$\pi(x) = 4c\min\{(x + \tfrac{2}{3})^\alpha, (\tfrac{4}{3} - x)^\alpha\} + c\min\{(x + \tfrac{1}{3})^\alpha, (\tfrac{5}{3} - x)^\alpha\}$$

for constants $c$ and $\alpha$ and for $0 < x < 1$. See Figure 1. We will compare three different proposal functions defined in $0 < x < 1$. Let the independence sampler have proposal function $q^1(x) = 1$. The random-walk sampler is given the proposal function $q_i^2(x) = 1/L$ for $x_{i-1} - L/2 < x < x_{i-1} + L/2$, and $q_i^2(x) = 0$ otherwise, where $x_{i-1}$ is the present position of the chain and $L$ is the maximum step length. In the case $x_{i-1}$ is close to 0 or 1, it is necessary to reduce $L$ in this iteration in order for the proposal function to be a proper density. This is very unlikely to happen except in the first few iterations. Before we can describe the proposal function for the adaptive independent sampler we need to define the two variables $z_i^j \in \tilde{y}^{i-1}$ such that $f(z_i^1) \geq f(z)$ for all $z \in \tilde{y}^{i-1}$ and $z_i^1, z < 0.5$ and that $f(z_i^2) \geq f(z)$ for all $z \in \tilde{y}^{i-1}$ and $z_i^2, z > 0.5$. Hence, the $z_i^j$'s are our best guess on the two modes after $i$ iterations of the Markov chain. The adaptive independent sampler is then given the proposal function

$$q^3(x)_i = \begin{cases} 1 - 2p + p/L, & \text{if } z_i^j - L/2 < x < z_i^j + L/2 \text{ for } j = 1,2, \\ 1 - 2p, & \text{otherwise,} \end{cases}$$

for two constants $p$ and $L$. $2p$ is the probability for a proposing a local jump and $L$ is the maximum length of a local jump. Also for this sampler it may be necessary to reduce $L$ in some iteration in order for the proposal function to be a proper density.

In the numerical calculations we set $\alpha = 2000$; see Figure 1. Then the modes are so sharp that 0.996 of the probability mass of $\pi(x)$ is located in two small intervals close to the modes with total length 0.01. Even though this example is schematic we believe it is quite representative for many MCMC problems. The typical MCMC is a random walk making small steps in each iteration. Often it is not critical whether this is made in one dimension at a time or there is only one dimension. We set the maximum local step length $L = 0.02$ in order to get an acceptance rate of about 0.25. If the step length is larger, the algorithm finds the mode faster but the acceptance rate decreases. The random-walk algorithm is very slow to move between the modes. This is illustrated in the lower left part of Figure 2 showing the minimum of the iteration number and the number of iterations since the chain was closer to the other mode. For random walk this is almost the same as the iteration number. Then the algorithm uses a very long time to find out that the probability mass close to the mode at $x = 1/3$ is larger. This slows down the convergence after about 100 iterations. The adaptive chain finds the modes faster than random walk. This is seen both in the convergence and in the ratio $\pi(1/3)/p_i(1/3)$ in the lower right part of the figure. Another important property is that the adaptive chain jumps easily between the modes. If we increase $L$, then the adaptive chain finds the mode faster, but does not jump as easily between the modes. We set $p = 0.4$ in order to jump often between the possible modes. The independence sampler uses many iterations to find the modes and the acceptance rate is very small. If we increase $\alpha$ making the modes sharper, then the difference between the different algorithms becomes even larger.

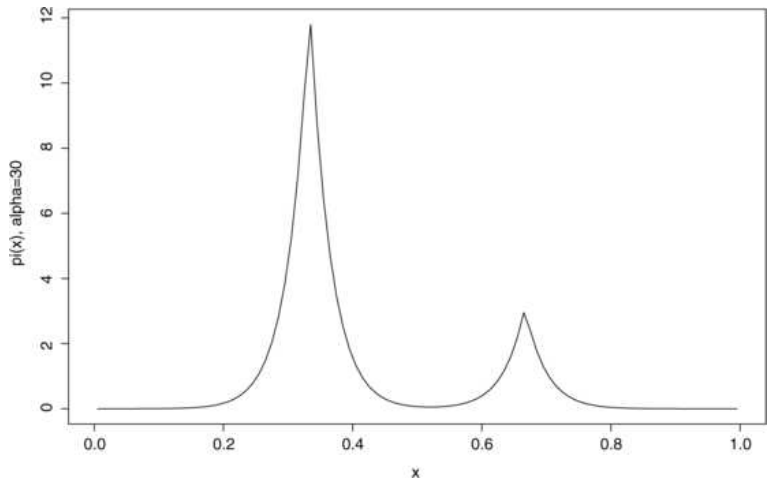

Fig. 1. *The function $\pi(x)$ in Example 1 with $\alpha = 30$.*

EXAMPLE 2. The adaptive chain is tested on a multimode example from Tjelmeland and Hegstad (2001). Let $\Omega = \mathbf{R}^n$ and $f(x) = \sum_{j=1}^{k} \omega_j \varphi_{\mu_j, \Sigma_j}(x)$ where $\varphi_{\mu_j, \Sigma_j}$ is the normal density in $\mathbf{R}^2$ with expectation $\mu_j$ and correlation matrix $\Sigma_j$. The constants $\omega_j > 0$ satisfy $\sum_j \omega_j = 1$. A natural proposal function to use in an adaptive Metropolis–Hastings algorithm, both here and generally, is a mixture of normal distributions:

$$(5) \qquad q_i^3(z|x, \tilde{y}^{i-1}) \propto \tau_0 \varphi_{\nu_0, \Lambda_0}(z) + \sum_{j=1}^{m_i} \tau_{i,j} \varphi_{\nu_{i,j}, \Lambda_j}(z).$$

The expectation $\nu_0$ is positioned central in the distribution. The corresponding variance $\Lambda_0$ is large. The other expectations $\nu_{i,j}$ for $j > 0$ are estimates of where the undersampling of $\pi$ is largest. The weights $\tau_{i,j}$ determine how often each distribution in the mixture is used for the proposal.

During the simulation we update a list of possible modes $\{\nu_{i,j}\}_{j=1}^{n_i} \subset \tilde{y}^{i-1}$. The list is empty when we start and is updated by the algorithm below based on the function $R(y) = f(y)/\varphi_{\nu_0, \Lambda_0}(y)$ with invariant $R(\nu_1) \geq R(\nu_2) \geq \cdots$. There is a maximum length of the list, that is, $n_i \leq M$. The most recent state in the history $y$ is considered included in the list by the algorithm:

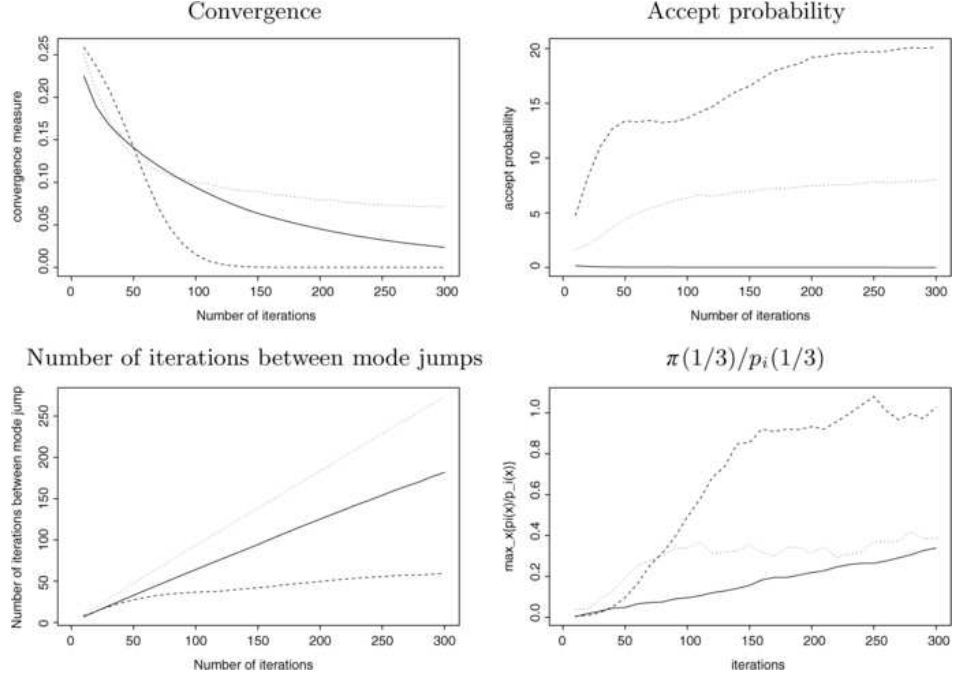

FIG. 2.  *The figure compares the three Metropolis–Hastings algorithms: independence sampler (line), random walk (dotted) and adaptive independent sampler (dashed). We see that the adaptive chain converges faster, jumps faster between the modes and obtains the correct density in the mode faster than the two other algorithms. Random walk has a higher acceptance rate, but with small steps each time. The simulation is based on 10,000 samples.*

1. If $R(y) > R(\nu_M)$ or $n_i < M$ then:
    (a) For $j = 1, \ldots, \min\{n_i, M\}$
        i. If $R(y) > R(\nu_j)$
            A. $<$Include $y$ in list before $\nu_j>$
            B. For $k = j+1, j+2, \ldots, \min\{n_i, M\}$

                If $\|y - \nu_k\| < \epsilon_1/2$ then

                    $<$Remove $\nu_k$ from list and exit loop$>$

        ii. Else if $\|y - \nu_j\| < \varepsilon_1$ exit loop.

The index $i$ is omitted in the algorithm and $\|\cdot\|$ denotes the Euclidean norm. Only $m_i \le M$ of the $\nu_j$ are used in the proposal function. More states are kept in the list since a new mode may remove several old ones. Let $m_i = \min\{M_0, n_i\}$ and $\tau_0 = 0.5$. The weights $\tau_{i,j}$ for $j = 1, 2, \ldots, m_i$ are defined as $\tau_{i,j} = 1/(5M_0) + cf(\nu_{i,j})$ where $c$ is defined such that $\sum_j \tau_{i,j} = 1$.

In addition, it makes sense to actively decrease the proposal probability in areas where previous proposals have shown that $f$ is small. This is especially attractive when evaluation of $f$ is expensive. A list of states $\xi_{i,j}$ are found, similarly as $\nu_{i,j}$ are found, maximizing $S(y) = \varphi^p_{\nu_0,\Lambda_0}(y)/f(y)$ for an exponent $p > 1$. The constants $N$, $N_0$ and $\varepsilon_2$ are defined as $M$, $M_0$ and $\varepsilon_1$, respectively. We can now actively reduce the sampling likelihood for the first term in (5) by introducing a new factor $\rho_i(z, \tilde{y}^{i-1})$, giving the proposal

$$(6) \qquad q_i^4(z|x, \tilde{y}^{i-1}) \propto \tau_0 \rho_i(z, \tilde{y}^{i-1}) \varphi_{\nu_0,\Lambda_0}(z) + \sum_{j=1}^{m_i} \tau_{i,j} \varphi_{\nu_{i,j},\Lambda_j}(z).$$

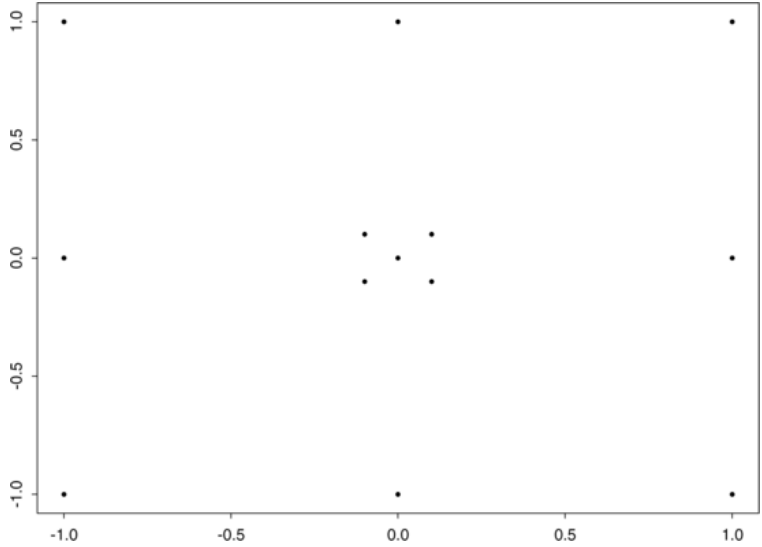

FIG. 3.  *Location of modes in the mixed normal distribution used as target distribution.*

The function $\rho_i(z, \tilde{y}^{i-1})$ is equal to a constant $c_i > 1$ except if the function $S(y)$ has been evaluated for $y$ close to $z$ and had a small value. In the latter case $\rho_i(z, \tilde{y}^{i-1})$ is small. Define the criteria $A_{i,t}$ so that $\|z - \xi_{i,t}\| < \varepsilon_2$ and $\|z - \xi_{i,k}\| \geq \varepsilon_2$ for all $k < t \leq N_0$. The following definition for $\rho$ is used:

$$\rho_i(z, \tilde{y}^{i-1}) = \begin{cases} \max\{\delta, f(\xi_{i,t})/\varphi_{\nu_0,\Lambda_0}(\xi_{i,t})\}, & A_{i,t} \text{ satisfied,} \\ c_i, & A_{i,t} \text{ not satisfied for any } t. \end{cases}$$

The constant $c_i$ is adjusted during the simulation in order to ensure that approximately 50% of the proposals come from the first term of the proposal function (6). Careful considerations on the smoothness of $f$ should be made when choosing $\epsilon_2$, since decreasing the proposal probability in an area also may worsen $\inf\{q_i(z|x, \tilde{y}^{i-1})/\pi(z)\}$, and hence worsen convergence.

The test example is the same as in Tjelmeland and Hegstad (2001) with $k = 13$ and modes $\mu_j$ located as shown in Figure 3. Each mode $i$ has the same weight $\omega_i = 1/13$ and covariance $\Sigma_i = \mathrm{diag}(0.01^2, 0.01^2)$. The variance is so small that a random walk between the outer modes is very unlikely, whereas random walks between the five center nodes are plausible. The constants in the adaptive chain are $M_0 = 20$, $M = 25$, and $\varepsilon_1 = 0.05$. The variances are set equal to $\Lambda_0 = \operatorname{diag}(1,1)$ and $\Lambda_j = \operatorname{diag}((0.03)^2, (0.03)^2)$ for $j > 0$. Two adaptive chains are evaluated: one with $\rho \equiv 1$ ($q^3$) and another ($q^4$) where $\rho$ depends on the $\xi_{...}$ list using $N_0 = N = 1000$, $\varepsilon_2 = 0.05$, $\delta = 0.1$ and $p = 1.3$. In addition two Metropolis–Hastings algorithms are implemented as a comparison with an independent proposer $q^1(y) = \varphi_{\nu_0, \Lambda_0}(y)$ and a random-walk proposer $q^2(y|x_{i-1}) = \varphi_{x_{i-1},(0.3)^2 \Lambda_0}(y)$, respectively.

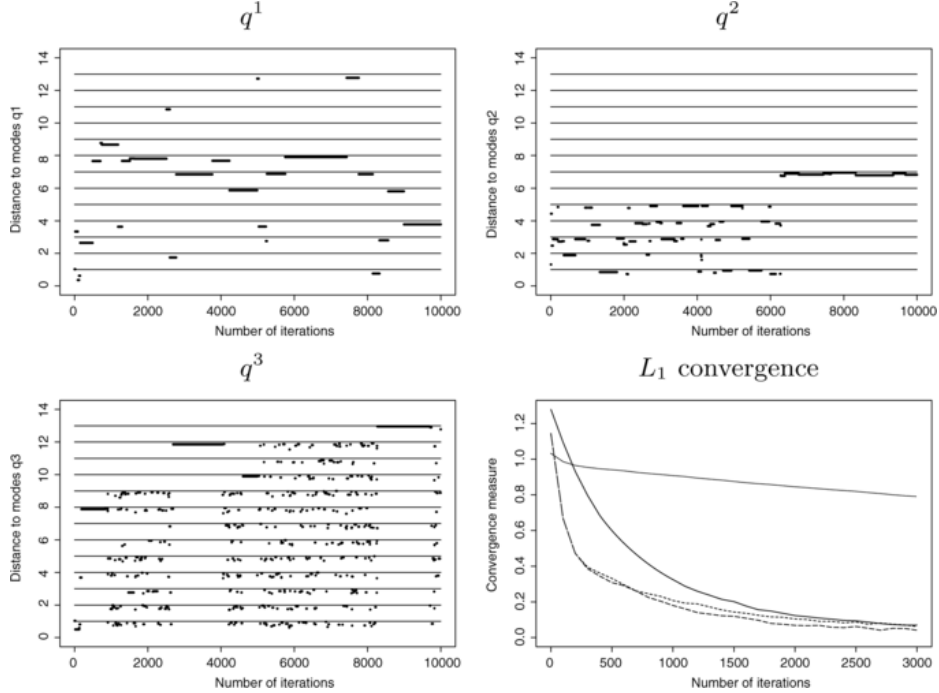

FIG. 4. *The chains: independent sampler ($q^1$, upper left), random walk ($q^2$, upper right), adaptive chain ($q^3$, lower left) jump between the 13 modes. The distance between a line and a dot indicates the distance between the sample and the mode. The modes closest to origin have smallest y value. The adaptive chain using the $\xi_{...}$ list converges fastest (lower right), closely followed by the other adaptive chain, then the independent sampler and finally the random-walk sampler. The last one is very slow since the probability of sampling far from origin is very small.*

Figure 4 shows the results from the experiments. The two adaptive chains are quite similar with acceptance rates equal to 0.07 and 0.09, respectively. The Markov chains have acceptance rates 0.001 and 0.004, respectively. The random walk jumps more often, but less often between modes. The adaptive chains resample often in modes that are identified earlier in the chain. The first time a mode is sampled, the chain leaves this mode with a small probability. Note also that the independent Markov chain does not sample as close to the mode value as the adaptive chain. This is the main reason for the difference in convergence, and is due to more sampling in the modes by the adaptive chain. The convergence is evaluated by calculating $\sum_j |r_j - 1/52|$ where $r_j$ is the ratio of chains within 52 regions which have the same probability in the target distribution $f$. There are four regions around each mode and each of these regions is described by the distance to the mode. The convergence measure using 20,000 chains will not be below 0.05 due to noise. The adaptive chain converges faster than the independent sampler. More important is the difference in mixing. In order to get 1000 samples with low correlation, the adaptive chain needs approximately 20,000 iterations (burn in 3000 and then sample every 10 iterations) and the independent Markov chain needs in the order of 1,000,000 (sample every 1000 iterations). There are several possibilities for further improvement. The variances $\Lambda_j$ and the constants $\varepsilon_k$ and $\delta$ may be adapted based on $\tilde{y}^{i-1}$. It is possible to run a local optimizer starting with some of the states in the $\nu_{..}$ list and include the endpoints of these optimizations into the list. The states in the $\xi_{..}$ list could be sorted into regions in order to avoid checking the complete list. Which variant of the algorithm is the best depends on the properties of $f$, in particular the cost of evaluating the function, in addition to the CPU resources available. The increased CPU resources required by changing from a traditional Metropolis–Hastings algorithm to an adaptive independent Metropolis–Hastings algorithm with $q^3$, are not large, but involves some programming. To use the proposal function $q^4$ increases the CPU resources considerably, and should only be used if the evaluation of $f$ is very expensive.

Note that even though the target distribution here was on the same form as the proposal, this method of adapting a proposal distribution is of general value. It can be viewed as a nonparametric approximation of the target distribution, bearing some resemblance to kernel methods.

EXAMPLE 3.  Let $f$ be the Cauchy distribution, $f(x) = 1/(1+x^2)\pi$. This distribution has heavier tails than the normal distribution and is known to give problems in the simulation; see Roberts and Stramer (2002). The same simulation algorithm as in Example 2 is used in this example. Only the adaptive chain $q^3$ and the independent Metropolis–Hastings algorithm $q^1(y) = \varphi_{0,1}(y)$ are tested. The constants in the adaptive chain are $M_0 = 70$, $M = 80$ and $\varepsilon_1 = 0.05$. The variances are set equal to $\Lambda_0 = 1$ and $\Lambda_j = (0.5)^2$ for $j > 0$. The modes $\nu_{i,j}$ in the proposal function will be in the tails of the distribution of $f$ with distance approximately $\varepsilon_1$ between neighboring modes. Figure 5 shows one chain for the two methods. Figure 6 shows the convergence estimated by dividing the state space into 20 equally likely regions for $|x|$ and uses the same norm as in the previous example. Ten thousand chains give a Monte Carlo noise level of approximately 0.07. The acceptance rates are 0.7 and 0.8, respectively, slightly higher for the adaptive chain. The adaptive chain reduces the problem with heavy tails by positioning modes in this area. This shows that the adaptive chain may be useful also when there are large differences between the limiting function and the proposal function. Both the adaptive and the independent chain converge even though the Doeblin condition is not satisfied. The independent sampler converges very slowly since the probability for sampling a state far from origin is too small. When such a state is sampled, the probability for leaving the state is very small.

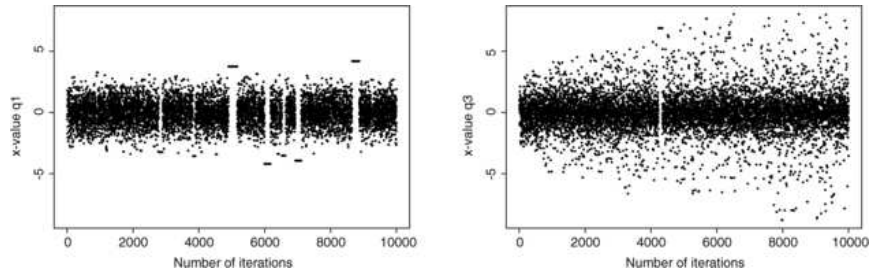

FIG. 5. *The Markov chain ($q^1$, left) samples the tails too seldom. The adaptive chain ($q^3$, right) increases the sampling of the tails as the modes $\nu_{...}$ identify these areas.*

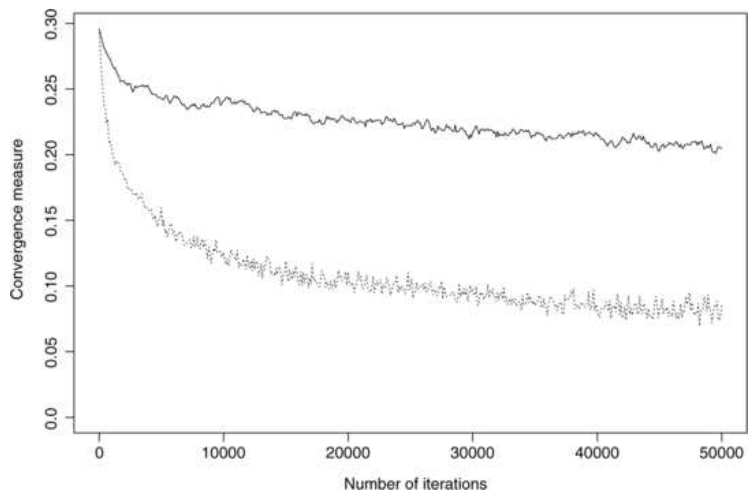

FIG. 6. *The adaptive chain converges much faster than the Markov chain. The adaptive chain reaches the Monte Carlo noise level for the convergence measure.*

EXAMPLE 4. In many applications in, for example, climate or petroleum, there are large simulation models. If each simulation takes several hours of CPU, it is necessary to carefully use all available information including previous runs before we start new simulations. There is often uncertainty in some input parameters and we want to find distributions for the input parameters that satisfy data that is output from the simulation model. It is outside the scope of this paper to describe such a large simulation model and we will illustrate this problem with a more schematic example.

We illustrate the simulation model with the function

$$f(\mathbf{x}) = 3\sin(x_1\pi) - x_1/2 + \sum_{i=2}^{5} \sin(x_i\pi/2)$$

assuming $0 < x_i < 1$ for $i = 1, 2, 3, 4, 5$. We use a Bayesian approach with a noninformative prior and the likelihood

$$l(\mathbf{x}|d) = \exp(-(f(\mathbf{x}) - d)^2/\sigma^2).$$

Then we want to simulate $\mathbf{x}$ proportional to $l(\mathbf{x}|d)$, but evaluate the function $f(\mathbf{x})$ as few times as possible. The function $f(\mathbf{x})$ may be approximated by the linear regression

$$\hat{f}(\mathbf{x}) = a_0 + \sum_{i=1}^{5} a_i x_i + b x_1^2$$

where the constants $a_i$ for $i = 0, 1, 2, 3, 4, 5$ and $b$ are evaluated from the evaluations of the function $f(\mathbf{x})$. Notice that most evaluations will be in the area where the likelihood is highest, that is, the linear regression will be best in the area of largest interest. In the adaptive independent sampler the proposal functions consist of two steps. We first propose uniformly $\mathbf{x} \in (0, 1)^5$ and then accept the proposal proportional to the function

$$\hat{l}(\mathbf{x}|d) = \exp(-(\hat{f}(\mathbf{x}) - d)^2/(5\sigma^2)).$$

Notice that we have included the factor 5 in the exponent in order to propose from a slightly larger area than the likelihood. This algorithm is compared with the independent sampler proposing uniformly from $\mathbf{x} \in (0, 1)^5$ and the random-walk sampler where each component of $\mathbf{x}$ is proposed changed uniformly in an interval with length $L = 0.1$ relative to the current state of the chain. The length of the interval is reduced close to the boundary of the domain.

We set $d = 2.5$ and $\sigma^2 = 0.005$. Then the likelihood is largest at two four-dimensional planes intersecting the domain close to two of the corners of the domain of $f$. The lower right figure in Figure 7 shows the likelihood along a line intersecting these two planes. The simulation shows that the adaptive sampler converges fastest and jumps easily between the two modes. The independent sampler jumps seldom and the random-walk sampler is not able to move between the modes and hence converge.

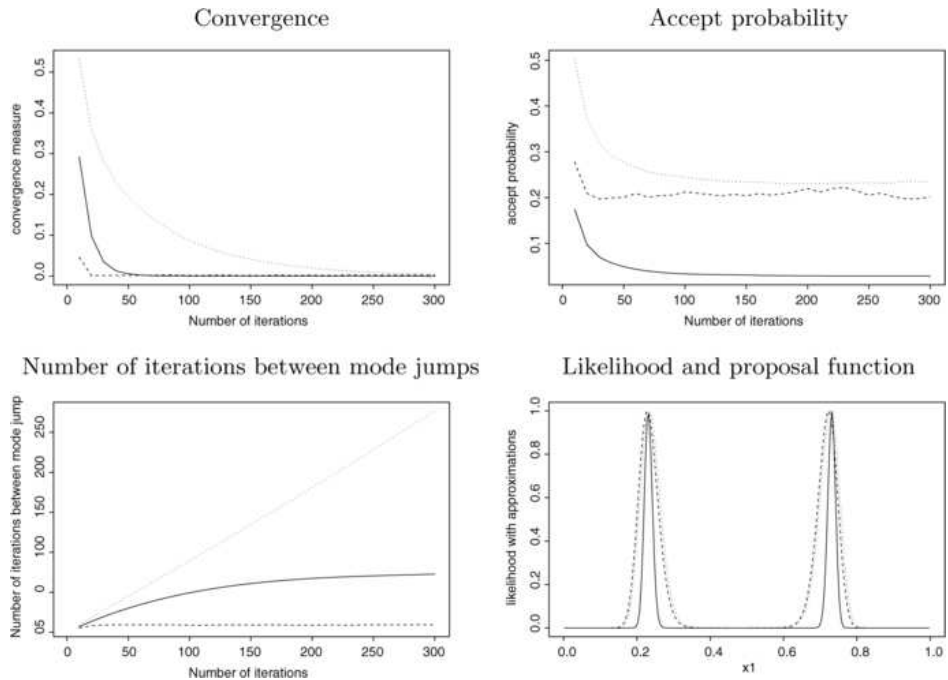

Fig. 7.  *The figure compares the three Metropolis–Hastings algorithms: independence sampler (line), random walk (dotted) and adaptive independent sampler (dashed). We see that the adaptive chain converges faster and jumps faster between the modes. Random walk has a higher acceptance rate, but with small steps each time. The simulation is based on 50,000 samples. The lower right figure shows the likelihood $l(x_1, 0.1, 0.1, 0.1, 0.1, 0.1|d)$ and the proposal function based on 10 and 50 evaluations of the simulation model along a line crossing the two planes with the modes.*

The traditional statistical approach for this problem would be to first perform some simulations of the large model estimating the linear regression parameters and then use the linear regression in the proposal function in a Metropolis–Hastings sampler. This implies to use no information in the data collection step and not the new information gained during the Metropolis–Hastings sampling. The adaptive sampler always uses all the information and all the time focuses on the most interesting area based on the available information. In a practical problem we will often not know the number of modes and whether we have found all the modes or not. The industrial practice in many such problems is often to neglect the uncertainty in the parameters and estimate the parameters using a standard gradient search optimization algorithm.

**5. Concluding remarks.**  The algorithm presented here allows fairly general adaption, based only on the assumption that a proposal distribution satisfying a strong Doeblin condition can be found. Convergence is geometric, with a rate that increases as the proposal distribution gets closer to the target distribution. The chain is also invariant for the limiting distribution in contrast to the more general adaptive chain. All previously proposed states, except the current state, can be used to generate the new proposal, allowing almost all previously gained information about the target distribution to be used. The algorithm is tested in four schematic examples, three with several modes and one with a heavy tail distribution, situations that are generally difficult for MCMC methods. Three examples are nonparametric and one is parametric. In all cases, the adaptive algorithm performed better than the standard MCMC alternatives we used for comparison.

The adaptive independent Metropolis–Hastings algorithm is a special case of the adaptive chain described in Holden (1998a), where convergence of the adaptive chain was proved by assuming that the detailed balance condition was satisfied for each possible history. This did in general not give invariance, which is a property of the adaptive independent Metropolis–Hastings. Convergence is here proved without assuming the detailed balance condition. Instead, we use rejection sampling and a restricted history which assures stationarity when the target distribution is reached.

The proposed algorithm may be considered as a generalization of the commonly used Metropolis–Hastings algorithm. If the proposal function in an iteration does not depend on the present state, then the present or the proposed state may be used for improving later proposal functions. Which state that may be used depends on whether the proposed state is accepted or not. This simple property contributes to the understanding of the commonly used Metropolis–Hastings algorithm. The paper also shows the close relationship between rejection sampling, MCMC and adaptive chains.

## APPENDIX: PROOFS

PROOF OF THEOREM 1. If the proposal function $q_{i+1}$ depends on the present state $x_i$, then the history is not extended and the theorem follows from standard Metropolis–Hastings theory. We will therefore focus on the case where $q_{i+1}$ is independent of the present state $x_i$, and we will use the notation $q_{i+1}(z_{i+1}|\tilde{y}^i)$.

Assume $\hat{p}_i(x_i|\tilde{y}^i) = \pi(x_i)$. Let $w$ be the state added to the history in step $i+1$, that is, $\tilde{y}^{i+1} = (\tilde{y}^i, w)$. Then $w \in \Omega$ is either $x_i$ or $z_{i+1}$ depending on whether the proposal $z_{i+1}$ is accepted or rejected in the $(i+1)$th iteration. The joint distribution of $(\tilde{y}^{i+1}, x_{i+1})$ is

$$\hat{p}_{i+1}(x_{i+1}|\tilde{y}^{i+1})\tilde{p}_{i+1}(\tilde{y}^{i+1})$$
$$= \tilde{p}_i(\tilde{y}^i)(\pi(w)q_{i+1}(x_{i+1}|\tilde{y}^i)\alpha_{i+1}(x_{i+1}, w, \tilde{y}^i)$$

$$+ \pi(x_{i+1})q_{i+1}(w|\tilde{y}^i)(1 - \alpha_{i+1}(w, x_{i+1}, \tilde{y}^i)))$$

$$= \tilde{p}_i(\tilde{y}^i)(\pi(x_{i+1})q_{i+1}(w|\tilde{y}^i)$$

$$+ \pi(w)q_{i+1}(x_{i+1}|\tilde{y}^i)\alpha_{i+1}(x_{i+1}, w, \tilde{y}^i)$$

$$- \pi(x_{i+1})q_{i+1}(w|\tilde{y}^i)\alpha_{i+1}(w, x_{i+1}, \tilde{y}^i))$$

$$= \pi(x_{i+1})\tilde{p}_i(\tilde{y}^i)q_{i+1}(w|\tilde{y}^i).$$

This shows that the chain never leaves the stationary distribution once it is reached and we have $\hat{p}_{i+1}(x_{i+1}|\tilde{y}^{i+1}) = \pi(x_{i+1})$.

Except for minor changes in the notation, the calculation above may also be used to prove the theorem when $q_{i+1}$ depends on the present state $x_i$. Then we must integrate over $w$ to get $\hat{p}_{i+1}(x_{i+1}|\tilde{y}^i) = \pi(x_{i+1})$. This is why the history cannot be extended in these iterations, and we set $\tilde{y}^{i+1} = \tilde{y}^i$. The critical point in the entire paper is the above calculation. Notice that it is essential that the proposal function does not depend on the present state, but may depend on the entire history $\tilde{y}^i$. $\quad\square$

PROOF OF THEOREM 2. Given that the limiting density is invariant for the algorithm conditioned on the history, the next step is to prove that a chain can reach the stationary distribution. This is done by observing that each iteration with an independent proposal also can be viewed as an iteration of a rejection sampler. Let $u_i$ be distributed uniformly between 0 and 1, and the proposal $z_i$ be accepted if $u_i < \alpha_i$. Let $A_i$ be the condition that $u_i q_i(z_i|\tilde{y}^{i-1})/\pi(z_i) \le a_i(\tilde{y}^{i-1})$. Then the distribution of $z_i$ given by $\tilde{y}^{i-1}, x_{i-1}$ and $A_i$ is proportional to

$$\text{Prob}(A_i|\tilde{y}^{i-1}, x_{i-1}, z_i)q_i(z_i|x_{i-1}, \tilde{y}^{i-1})$$

$$= \text{Prob}\left(u_i \le a_i(\tilde{y}^{i-1})\frac{\pi(z_i)}{q_i(z_i|x_{i-1}, \tilde{y}^{i-1})}\right)q_i(z_i|x_{i-1}, \tilde{y}^{i-1})$$

$$= a_i(\tilde{y}^{i-1})\pi(z_i).$$

This means that if $A_i$ is satisfied, $z_i$ is a sample from $\pi$. Furthermore, if $A_i$ is satisfied, $u_i \le \alpha_i$, and hence $z_i$ is always accepted. If $q_i$ does not depend on $x_{i-1}$, then $x_{i-1}$ is appended to the history, and we conclude that $\hat{p}_i(x|\tilde{y}^i, A_i) = \pi(x)$. If $q_i$ depends on $x_{i-1}$, we must integrate over the distribution of $x_{i-1}$ given $\tilde{y}^i, A_i$ to obtain the same conclusion. Finally, $A_i$ is satisfied with probability $a_i(\tilde{y}^{i-1})$. This implies that the probability of jumping to the limiting distribution in iteration $i$ is $a_i(\tilde{y}^{i-1})$ independent of the distribution of $x_i$. From Theorem 1 we see that when we reach the limiting density we remain in the limiting density. Hence, this is an absorbing state with a probability larger than $a > 0$ to jump to this in each iteration. Then the probability to be in the limiting density within a finite number of samples may be made arbitrarily close to 1.

Let

$$b_i(\tilde{y}^{i-1}) = \prod_{j=1}^{i}(1 - a_j(\tilde{y}^{j-1})).$$

The probability for $x_i$ being from the limiting distribution is then described by the following chain: Let $I_0 = 0$ and for $i \geq 0$

$$I_{i+1} = \begin{cases} 0, & \text{with probability } 1 - a_{i+1}(\tilde{y}^i) \text{ if } I_i = 0, \\ 1, & \text{otherwise.} \end{cases}$$

Clearly

$$\text{Prob}(I_i = 0|\tilde{y}^{i-1}) = \prod_{j=1}^{i}(1 - a_j(\tilde{y}^{j-1})) = b_i(\tilde{y}^{i-1}).$$

This implies that

(7) $$\hat{p}_i(x_i|\tilde{y}^i) = \pi(x_i)(1 - b_i(\tilde{y}^{i-1})) + v_i(x_i|\tilde{y}^i)b_i(\tilde{y}^i)$$

where $v_i$ is a distribution. This gives the following bound on the error:

$$\|p_i - \pi\|_{\text{TV}} = \int_\Omega \left| \int_{\Omega^i} \hat{p}_i(x_i|\tilde{y}^i)\tilde{p}_i(\tilde{y}^i)\,d\mu(\tilde{y}^i) - \pi(x_i) \right| d\mu(x_i)$$

$$= \int_\Omega \left| \int_{\Omega^i} (\pi(x_i)(1 - b_i(\tilde{y}^{i-1})) + v_i(x_i|\tilde{y}^i)b_i(\tilde{y}^{i-1}) - \pi(x_i)) \right.$$

$$\left. \times \tilde{p}_i(\tilde{y}^i)\,d\mu(\tilde{y}^i) \right| d\mu(x_i)$$

$$= \int_\Omega \left| \int_{\Omega^i} (v_i(x_i|\tilde{y}^i) - \pi(x_i))\tilde{p}_i(\tilde{y}^i)b_i(\tilde{y}^{i-1})\,d\mu(\tilde{y}^i) \right| d\mu(x_i)$$

$$\leq \int_{\Omega^i} \int_\Omega |v_i(x_i|\tilde{y}^i) - \pi(x_i)|\,d\mu(x_i)\,\tilde{p}_i(\tilde{y}^i)b_i(\tilde{y}^{i-1})\,d\mu(\tilde{y}^i)$$

$$\leq 2 \int_{\Omega^i} \tilde{p}_i(\tilde{y}^i)b_i(\tilde{y}^{i-1})\,d\mu(\tilde{y}^i),$$

proving (3).

If $a_{j+1}(\tilde{y}^j) \geq a_{j+1}$ for all $j$ and $\tilde{y}^j \in \Omega^j$, then (4) follows trivially. Note that the probability for jumping to the limiting distribution in iteration $i$ is $a_i(\tilde{y}^{i-1})$ independent of the distribution of $x_i$. $\quad\square$

**Acknowledgments.** The authors thank Arnoldo Frigessi, Håkon Tjelmeland and Øivind Skare for valuable comments and contributions.

# REFERENCES

ANDRIEU, C. and MOULINES, E. (2003). On the ergodicity properties of some adaptive Monte Carlo algorithms. Technical report, Univ. Bristol.

ATCHADE, Y. and ROSENTHAL, J. (2003). On adaptive Markov chain Monte Carlo algorithms. Technical report, Univ. Ottawa.

ERLAND, S. (2003). *Approximating hidden Gaussian Markov random fields.* Ph.D. thesis, Dept. Mathematical Sciences, Norwegian Univ. Science and Technology, Trondheim, Norway.

GÅSEMYR, J. (2003). On an adaptive version of the Metropolis–Hastings algorithm with independent proposal distribution. *Scand. J. Statist.* **30** 159–173. MR1965100

GEYER, C. J. (1992). Practical Markov chain Monte Carlo. *Statist. Sci.* **7** 473–483.

GILKS, W. R., RICHARDSON, S. and SPIEGELHALTER, D. J., eds. (1996). *Markov Chain Monte Carlo in Practice.* Chapman & Hall, London. MR1397966

GILKS, W. R., ROBERTS, G. O. and SAHU, S. K. (1998). Adaptive Markov chain Monte Carlo through regeneration. *J. Amer. Statist. Assoc.* **93** 1045–1054. MR1649199

HAARIO, H., SAKSMAN, E. and TAMMINEN, J. (2001). An adaptive Metropolis algorithm. *Bernoulli* **7** 223–242. MR1828504

HOLDEN, L. (1998a). Adaptive chains. NR-note SAND/11/1998, Norwegian Computing Center. Available at http://www.statslab.cam.ac.uk/~mcmc/.

HOLDEN, L. (1998b). Geometric convergence of the Metropolis–Hastings simulation algorithm. *Statist. Probab. Lett.* **39** 371–377. MR1646224

KLEIJNEN, J. P. C. and SARGENT, R. G. (2000). A methodology for fitting and validating metamodels in simulation. *European J. Oper. Res.* **120** 14–29.

LIU, J. S. (1996). Metropolized independent sampling with comparison to rejection sampling and importance sampling. *Statist. Comput.* **6** 113–119.

MEYN, S. P. and TWEEDIE, R. L. (1993). *Markov Chains and Stochastic Stability.* Springer, London. MR1287609

PARK, J. S. and JEON, J. (2002). Estimation of input parameters in complex simulatioin using a gaussian process metamodel. *Probabilistic Engineering Mechanics* **17** 219–225.

ROBERTS, G. O. and ROSENTHAL, J. S. (1998). Markov-chain Monte Carlo: Some practical implications of theoretical results. *Canad. J. Statist.* **26** 5–31. With discussion by Hemant Ishwaran and Neal Madras and a rejoinder by the authors. MR1624414

ROBERTS, G. O. and ROSENTHAL, J. (2005). Coupling and ergodicity of adaptive MCMC. Technical Report 314, Univ. Toronto.

ROBERTS, G. O. and STRAMER, O. (2002). Tempered Langevin diffusions and algorithms. Technical Report 314, Univ. Iowa.

TIERNEY, L. (1994). Markov chains for exploring posterior distributions. *Ann. Statist.* **22** 1701–1762. With discussion and a rejoinder by the author. MR1329166

TJELMELAND, H. and HEGSTAD, B. K. (2001). Mode jumping proposals in MCMC. *Scand. J. Statist.* **28** 205–223. MR1844357

L. HOLDEN
R. HAUGE
M. HOLDEN
NORWEGIAN COMPUTING CENTER
P.O. BOX 114 BLINDERN
N-0314 OSLO
NORWAY
E-MAIL: Lars.Holden@nr.no
        Ragnar.Hauge@nr.no
        Marit.Holden@nr.no